\documentclass[a4paper,12pt]{article}
\usepackage[all]{xy}
\usepackage{amsmath, amsthm}
\usepackage{amssymb,amsfonts}

\date{}
\title{Admissibility of groups over function fields of p-adic curves}
\author{B. Surendranath Reddy and V. Suresh}

\begin{document}
\maketitle
\section*{Abstract}
Let $K$ be a field and $G$ a finite group. The question of `admissibility' of $G$
over $K$ was originally posed by Schacher,
who gave partial results in the case $K=\mathbb{Q}$. In this paper,
we give  necessary conditions for admissibility of a finite group
$G$ over function fields of curves over complete discretely
valued fields.  Using this criterion, we give an example of a finite group
which is not admissible over $\mathbb{Q}_{p}(t)$. We also prove a
certain Hasse principle for division algebras over such fields. 

\section*{Introduction }
Let $K$ be a field and $G$ a finite group. We say that $G$ is {\it
admissible} over $K$ if there exists a division ring $D$ central
over $K$ and a maximal subfield  $L$ of $D$ which is Galois over $K$
with Galois group $G$. Schacher asked, given a field $K$, which
finite groups are admissible over $K$ and proved that if a finite
group $G$ is admissible over $\mathbb{Q}$, then every Sylow
subgroup of $G$ is metacyclic ([Sc], 4.1). This led to the conjecture that
a finite group $G$ is admissible over $\mathbb{Q}$ if and only if
every Sylow subgroup of $G$ is meta-cyclic. This conjecture has
been proved for all solvable groups ([So]) and for certain
non-solvable groups of small order ([CS], [FS], [Fe1], [Fe2], [Fe3]).  

Recently Harbater, Hartman and Krashen ([HHK2], 4.5) gave a
characterization of admissible groups over function fields of curves
over complete discretely valued fields with algebraically closed
residue fields.   In
this paper, we consider the function fields of curves over complete
discretely valued fields without any assumptions on the residue
fields and prove the following

\paragraph*{Theorem 1.} Let $K$ be a complete discretely valued field
with residue field $k$ and $F=K(X)$ be the function field of a curve
$X$ over $K$. Let $G$ be a finite group. Suppose  that the order of $G$ 
is coprime to char$(k)$. If $G$ is admissible over $F$ then
every Sylow subgroup $P$ of $G$ has a normal series $P\supseteq
P_{1}\supseteq P_2$ such that
 
(1) $P/P_{1}$ and $P_{2}$ are cyclic

(2) $P_{1}/P_{2}$ is admissible over a finite extension of the residue
field of a discrete valuation of $F$.

\vskip 3mm

A main ingredient for the proof of the above theorem is  the following Hasse  
principle for central simple  algebras, which has independent interest. 

\paragraph*{Theorem 2.}  Let $K$ be a complete discretely valued field
with residue field $k$
and $F=K(X)$ be the function field of a curve $X$ over $K$. Let $A$
be a central simple algebra over $F$ of degree  $n = \ell^r$ for
some prime $\ell$ and $r \geq 1$. Assume that $\ell$ is not equal to 
 char$(k)$ and $K$ contains a primitive $n^{th}$ root of unity.
Then index$(A) = $ index$(A \otimes F_v)$   for some discrete valuation $v$ of
$F$.

\vskip 3mm

For the proof of the above theorem, we use  the patching techniques 
of ([HHK1]). A similar Hasse principle is proved for quadratic forms  over such  
fields in ([CTPS], 3.1). In ([HHK3], 9.12), it is proved that if a central simple algebra $A$ over
$F$ ($F$ as above), is split over $F_{\nu}$ for all discrete valuations on $F$,
then $A$ is split over $F$. 

There are some examples of classes of finite groups which are admissible
over the rational function fields. However there was no example, in
the literature,  of a
finite group which is not admissible over $\mathbb{Q}_p(t)$. 
Using  Theorem 1, we give an example of a finite group which is not admissible over 
$\mathbb{Q}_p(t)$. We also prove admissibility of a certain class of groups over
$\mathbb{Q}_{p}(t)$ using patching techniques.

\paragraph*{Theorem 3.} Let K be a $p$-adic field and $F$ the 
function field of a curve over $F$. 
Let $G$ be a finite group with order coprime to $char(k)$. If every
Sylow subgroup of $G$ is a quotient of $\mathbb{Z}^4$, then $G$ is
admissible over $F$.

\vskip 3mm

 In [FSS], it was proved that every abelian group on three or less
 generators is admissible over $\mathbb{Q}(t)$. We conclude by showing
 that every abelian group of order $n$  with four or less generators is admissible
 over $\mathbb{Q}(\zeta)(t)$, where $\zeta$ is a primitive $n^{th}$
 root of unity. 

\section*{1. Some Preliminaries}

In this section we recall a few basic definitions and facts about
division algebras and patching techniques ([GS], [HH], [HHK1], [P],
[S2],  [Sc], [Sch], [Ser]).

Let $K$ be a field and $Br(K)$ be the Brauer group of central simple
algebras over $K$. For an integer $n \geq 2$, let $_nBr(K)$ denote the
$n$-torsion subgroup of $Br(K)$. If $A$ and $B$ are two central simple
algebras over $K$, we write $A \simeq B$ if $A$ and $B$ are isomorphic
as $K$-algebras and we write $A = B$ if they represent the same
element in $Br(K)$. 
Let  $n$ be an integer coprime to char$(K)$.
Suppose  $E/K$ is a cyclic extension of degree $n$ and $\sigma$ a generator of
$Gal(E/K)$. For $b \in K^*$,  let $(E/K, \sigma, b)$ (or simply $(E,
\sigma, b)$) be the $K$-algebra
generated by $E$ and $y$ with $y^n = b$ and $\lambda y = y
\sigma(\lambda)$ for all $\lambda \in E$. Then $(E, \sigma, b)$ is a
central simple algebra over $K$ and represents an element in
$_nBr(K)$.  Suppose that $K$ contains a primitive $n^{th}$ root of
unity.  Then $E = K(\sqrt[n]{a})$ for some $a \in K^*$  and
$\sigma(\sqrt[n]{a}) = \zeta \sqrt[n]{a}$ for a
primitive $n^{th}$ root of unity $\zeta \in K^*$.
The algebra  $(E, \sigma, b)$ is  generated by $x, y$
with $x^n = a$, $y^n = b$ and $xy = \zeta yx$ and we denote the
algebra  $(E, \sigma, b)$ by  $(a, b)_n$.

Suppose that $\nu$ is a discrete
valuation of $K$ with residue field $k$. Let $n$ be a natural
number which is coprime to  char$(k)$. 
Then we
have a {\it residue homomorphism } $\partial_{\nu} :~ _nBr(K) \to
H^1(k, \mathbb{Z}/n\mathbb{Z})$, where $H^1(k,
\mathbb{Z}/n\mathbb{Z})$ denotes the first Galois cohomology group.
Suppose $k$ contains a
primitive $n^{th}$ root of unity.  By fixing a primitive $n^{th}$
root of unity, we identify $H^1(k, \mathbb{Z}/n\mathbb{Z})$
with $k^*/k^{*n}$. With this identification we have 
$\partial_{\nu}((a, b)_n) = \overline
{\frac{a^{\nu(b)}}{b^{\nu(a)}}} \in k^*/k^{*n}$, where for
any $c \in K^*$ which is a unit at $\nu$, $\overline{c}$ denotes its
image in $k^*$. More generally, let $E/K$ be a cyclic unramified
inert extension of $K$ and    
$\sigma \in Gal(E/K)$ be a generator.  Let $\pi \in K^*$ be a parameter at $\nu$.  
Then the residue of $(E/K, \sigma,\pi)$ is
$(E_0/k, \sigma_0)$, where $k$ is the residue field of  $K$ 
at $\nu$,  $E_0$ is the residue field of $E$ at the unique extension
of $\nu$ to $E$ and $\sigma_0$ is the image of $\sigma$.

Let $K$ be a field and $n$ an integer. Then
$H^{1}(K,\mathbb{Z}/n\mathbb{Z})$ classifies the equivalence classes of pairs $(E, \sigma)$,
where $E$ is a cyclic Galois field
extension of $K$ of degree a divisor  of $n$ and $\sigma$ a generator
of the Galois group $G(E/K)$ of $E/K$. 
Let $(E,\sigma)$ be a pair  representing an element in $H^1(K, \mathbb{Z}/n\mathbb{Z})$.
Let $m \geq 1$ be an integer. We now describe $m(E, \sigma)$
in the group $H^1(K, \mathbb{Z}/n\mathbb{Z})$.  Let $d$ be the greatest common divisor
of $m$ and $[E: K]$. Let $E(m)$ be the subfield of $E$ fixed by
$\sigma^{[E:K]/d}$
and $\sigma(m)$ be the restriction of 
$\sigma^{m/d}$ to $E(m)$. Then $m(E, \sigma) = (E(m), \sigma^{m/d})$. 
The order of $(E, \sigma)$ in $H^1(K, \mathbb{Z}/n\mathbb{Z})$ is
equal to $[E : K]$. 
Since the identity element of $H^1(K, \mathbb{Z}/n\mathbb{Z})$ is $(K, id)$, we have 
$m(E, \sigma)$ is trivial in $H^1(K,\mathbb{Z}/n\mathbb{Z})$ 
if and only if $E(m) = K$ if and only if  $[E : K]$   divides $m$.

Let $K$ be a complete discretely valued field with residue field $k$.
Let $L/K$ be a finite extension of $K$. Since $K$ is complete, the
discrete valuation of $K$ extends uniquely to a discrete valuation of
$L$ and $L$ is complete with respect to this discrete valuation. 
Let $n$ be a natural number which is coprime to the characteristic
of $k$.   Let $L/K$ be a Galois  extension of degree $n$.  Let $L_1$
be the maximal unramified extension of $K$. Since $n$ is coprime to
char$(k)$, the residue field $L_0$ of $L$ is same as the residue field
of $L_1$. Since $L/K$ is Galois, $L_1/K$ and $L_0/k$ are also Galois
and there is a natural isomorphism $Gal(L/K) \to Gal(L_0/k)$.  
Let $E_0/k$ be a cyclic extension of degree $n$. Then there
exists a unique (up to isomorphism) unramified cyclic extension $E$ of $K$ with
residue field $E_0$ and a natural isomorphism $Gal(E/K) \to
Gal(E_0/k)$. Let $(E_0, \sigma_0) \in H^1(k,
\mathbb{Z}/n\mathbb{Z})$. Then we have a unique $(E, \sigma) \in
H^1(K, \mathbb{Z}/n\mathbb{Z})$ with $[E: K ] = [E_0 : k]$, 
$E_0$ as the residue field of $E$ and 
the image of $\sigma$ equal to $\sigma_0$; $(E, \sigma)$  is called the {\it lift} of $(E_0,
\sigma_0)$.

Let ${\cal X}$ be a regular integral scheme with function field $F$.
Let $n$ be an integer which is a unit on ${\cal X}$. Let $f \in F$
and $P \in {\cal X}$ be a point. If $f$ is regular at $P$, then we
denote its image in the residue field $\kappa(P)$ at $P$ by $f(P)$.
Let ${\cal X}^1$ denote the set of codimension one points of ${\cal
X}$. For each codimension one point $x$ of ${\cal X}$, we have
discrete valuation $\nu_x$ on $F$. Let $\kappa(x)$ denote the
residue field at $x$. Since $n$ is a unit on ${\cal X}$, $n$ is
coprime to char$(\kappa(x))$ and we have the residue
homomorphism $\partial_x : ~_nBr(F) \to H^1(\kappa(x),
\mathbb{Z}/n\mathbb{Z})$. Let $\alpha \in ~_nBr(F)$. We say that
$\alpha$ is {\it unramified} at $x$ if $\partial_x(\alpha) = 0$. We
say that $\alpha$ is {\it unramified  } on ${\cal X}$ if it is
unramified at every codimension one point of ${\cal X}$. Let $A$ be a
central simple algebra over $F$. We say tat $A$ is {\it unramified}
if its class in $Br(F)$ is unramified. If ${\cal X} = Spec(B)$,
then we say that $\alpha$ is unramified on $B$ if it is
unramified on ${\cal X}$.

Let $F$ be a field and $A$ a central simple algebra over $F$. 
Then $A \simeq M_m(D)$ for some central division algebra $D$ over $F$.
The {\it degree} of $A$ is defined as $\sqrt{dim_FA}$ and the {\it index} of $A$
is defined as the degree of $D$.  If $L/K$ is a field extension, then
 index$(A \otimes_FL)$ divides index$(A)$. Let $B$ be an integral domain and
$F$ its  field of fractions. Let ${\cal A}$ be an Azumaya algebra over
$B$, then we define the   {\it index} of ${\cal A}$ to be the index of ${\cal A}\otimes _B F$.

Let  $B$ be a regular integral domain of
dimension at most 2  which is complete with respect to a prime ideal
$P$.  Assume that $B/P$ is a regular integral domain of dimension at
most 1 (for example $P$ is  a  maximal ideal). 
Let $\kappa(P)$ be the field of fractions of $B/P$.  
Let $A$ be a central simple algebra over the field of fractions $F$ of
$B$  which is unramified on $B$. Then there exists an
Azumaya algebra ${\cal A}$ over $B$ such that ${\cal A}
\otimes_B F \simeq A$ ([CTS], 6.13)).   For an ideal $I$ of $B$, we denote the
algebra ${\cal A} \otimes_B B/I$ by $A(I)$.

\paragraph*{Lemma 1.1}  Let $A$, $B$, ${\cal A}$ and $P$ be as
above. Then index$(A) = $  index $({\cal
  A}  \otimes_{B/P} \kappa(P))$.

\paragraph*{Proof.}  Suppose $A \simeq M_n(D)$ for some division
algebra $D$ over $F$.  
Since $A$ is unramified  on $B$, $D$ is also unramified on $B$.
Since $B$ is a regular domain of dimension  2, there exists an
Azumaya algebra ${\cal D}$ over $B$ such that ${\cal D} \otimes_B K
\simeq D$ ([CTS], 6.13)).  Since $Br(B) \to Br(F)$ is
injective ([AG], 7.2), ${\cal A} = {\cal D}$ in  $Br(B)$.   In
particular $ {\cal A}
\otimes_{B/P} \kappa(P)  =  {\cal D} \otimes _{B/P}  \kappa(P)$  in
$Br(\kappa(P))$.
Hence index$(A) = $ degree$(D) = \sqrt{{\rm dim}_F(D)} = \sqrt{{\rm
    rank}_B({\cal D})} =\sqrt {{\rm rank}_{B/P}({\cal D} \otimes _B
  B/P)} = \sqrt{{\rm dim}_{\kappa(P)}({\cal D}  \otimes_{B/P} \kappa(P))}
\geq$ index$({\cal A} \otimes_{B/P} \kappa(P))$. 

Suppose ${\cal A} \otimes _{B/P} \kappa(P) \simeq M_m(D_0)$ for some
central division algebra $D_0$ over $\kappa(P)$.  Since $B/P$ is
regular integral domain of dimension  1, there exists an Azumaya
algebra ${\cal D}_0$ over $B/P$ such that ${\cal D}_0 \otimes _{B/P}
\kappa(P) \simeq D_0$.  Since $B$ is $P$-adically complete, there
exists an Azumaya algebra $\tilde{{\cal D}}_0$ over $B$ such that $\tilde{{\cal
  D}}_0\otimes _B  B/P \simeq {\cal D}_0$   and  $ {\cal A}  = \tilde{{\cal D}}_0$ in
$Br(B)$  ([C], [KOS]).  In particular $A   = \tilde{{\cal D}}_0
\otimes_B F$ in $Br(F)$.  Hence 
index$({\cal A} \otimes_{B/P} \kappa(P)) = $ degree$(D_0) = \sqrt{{\rm
    dim}_{\kappa(P)}(D_0)} =  \sqrt{{\rm
    rank}_{B/P}({\cal D}_0) } =\sqrt {{\rm rank}_B(\tilde{{\cal D}}_0)}
  =  \sqrt{{\rm dim}_F({\tilde{{\cal D}}_0}   \otimes  F)}
\geq$ index$(  A )$.  Thus index$(A) = $ index$({\cal A} \otimes_{B/P} \kappa(P))$.
\hfill $\Box$

\vskip 3mm

Let $K$ be a field and $L$ a finite extension of $K$. Then $L$ is
called $K$-\emph{adequate} if there is a division ring $D$ central
over $K$ containing $L$ as a maximal subfield. A finite group $G$ is
called $K$-\emph{admissible} if there is a Galois extension $L$ of
$K$ with $G$ as the Galois group of $L$ over $K$, and $L$ is
$K$-\emph{adequate}.

A finite group $G$ is called \emph{metacyclic} if $G$ has a normal
subgroup $H$ such that $H$ is cyclic and $G/H$ is cyclic.

\section*{2. A local-global principle for central simple  algebras}

Let $K$ be a complete discretely valued field with residue field $k$.
Let $F$ be the  function field of a curve over $K$. Let $n$ be an
integer which is coprime to the characteristic of $k$. Assume that
$K$ contains a primitive $n^{th}$ root of unity. In this section we
prove a certain Hasse principle for central simple algebras over $F$
of index $n$. We begin with the following

\paragraph*{Lemma 2.1.} (cf. [FS], Proposition 1(3) and [JW], 5.15)
Let $R$ be a complete discrete valuated ring and $K$ its field of
fractions. Let $A$ be a central simple algebra over $K$ of index
$n$ which is unramified at $R$. Let $E$ be an unramified  cyclic extension of $K$ of degree $m$
and $\sigma$ a generator of the Galois group of $E/K$. Let $\pi$ be
a parameter in $R$. Assume that $mn$ is invertible in
$R$. Then index$(A \otimes (E, \sigma, \pi))=$ index$(A \otimes E ) \cdot [E: K]$.
 
\vskip 3mm

The following two lemmas (2.2, 2.3) are well known.

\paragraph*{Lemma 2.2.} Let $R$ be a regular ring of
dimension 2 with field of fractions $F$. Let $n$ be an integer which
is  a unit in $R$. Assume that $F$ contains a primitive $n^{\rm th}$ root of unity.
Suppose $A$ is a central  simple algebra over $F$ which is unramified on $R$. 
Let $x \in R$ be a regular prime and $\nu$ be the  discrete valuation on $F$ given by $x$.
Suppose that $R$ is complete with respect to $(x)$-adic topology. Let $u \in R$ be a unit.
Then index$(A\otimes F_{\nu}(\sqrt[n]{u})) = $ index$(A \otimes F(\sqrt[n]{u}))$.

\paragraph*{Proof.} Let $S$ be the integral closure of $R$ in $F(\sqrt[n]{u})$. 
Since $R$ is a  regular  ring and $n$, $u$  are units in $R$, 
$S$ is also a regular  ring.  Let $x$ be a regular prime in 
$R$, then $x$ is also a regular prime in $S$. Thus by replacing $R$ by $S$, it is 
enough to show that  index$(A \otimes F_{\nu}) = $ index$(A)$.

Since $R$ is a two-dimensional  regular  ring, there is
an Azumaya algebra ${\cal A}$ over $R$ with ${\cal A} \otimes F \simeq A$ ([CTS], 6.13)).
Since $R$ is  complete with respect to $(x)$-adic topology, 
index$(A) =$ index$({\cal A}\otimes \kappa(x))$ (cf. 1.1), where
$\kappa(x)$ is the field of fractions of $R/(x)$.      
Let $R_{\nu}$ be the ring of integers in $F_{\nu}$. Since $A$ is unramified 
on $R$, $A$ is also unramified on $R_{\nu}$. Since $F_{\nu}$ is complete,  
index$(A \otimes F_{\nu})  = $ index$({\cal A}\otimes \kappa(\nu))$
(cf. 1.1). Since $\kappa(\nu) = \kappa(x)$,  
index$(A) =$ index$( A  \otimes F_{\nu} )$. 
\hfill $\Box$

\paragraph*{Lemma 2.3.}  Let $R$ be a complete regular local ring 
with field of fractions $F$ and residue field $k$. Let $n$ be an integer which
is  a unit in $R$. Let $u_1, \cdots , u_r \in R$ be  units. Suppose $x \in R$ is a regular prime.
Let $\nu$  be the  discrete valuation on $F$ given by $x$. Then 
$[F_{\nu}(\sqrt[n]{u_1}, \cdots , \sqrt[n]{u_r} ) : F_{\nu}] = [F(\sqrt[n]{u_1}, \cdots, \sqrt[n]{u_r}) : F]$.

\paragraph*{Proof.} Since $F_{\nu}$ is a complete discretely  valued field and $n$, $u_1$, $\cdots$, $u_r$  are units
in the ring of integers, $[F_{\nu}(\sqrt[n]{u_1}, \cdots , \sqrt[n]{u_r}) : F_{\nu}] = 
[\kappa(\nu)(\sqrt[n]{\overline{u}_1}, \cdots ,\sqrt[n]{\overline{u}_r}): \kappa(\nu)]$.
Since $\kappa(\nu)$ is the field of fractions of the complete local ring $R/(x)$ and the the residue
field  of $R/(x)$ is $k$,  we have $[\kappa(\nu)(\sqrt[n]{\overline{u}_1}, \cdots ,\sqrt[n]{\overline{u}_r}): \kappa(\nu)] = 
[k(\sqrt[n]{\overline{u_1}}, \cdots, \sqrt[n]{\overline{u_r}}) : k]$.  
Since $R$ is complete and $u_1, \cdots, u_r$ are units, we also have $[F(\sqrt[n]{u_1}, \cdots , \sqrt[n]{u_r}) : F]
= [k(\sqrt[n]{\overline{u}_1}, \cdots , \sqrt[n]{\overline{u}_r}) : k]$. 
Hence $[F_{\nu}(\sqrt[n]{u_1}, \cdots , \sqrt[n]{u_r}) : F_{\nu}] = [F(\sqrt[n]{u_1}, \cdots, \sqrt[n]{u_r}) : F]$.
\hfill $\Box$

\paragraph*{Proposition 2.4.} Let $R$ be a 2-dimensional complete regular local ring 
with maximal ideal $m = (x, y)$. Let $F$ be the field of fraction of $R$ and $k$ 
the residue field of $R$. Let $n$ be an integer coprime to char$(k)$. Assume that $F$ contains 
a primitive $n^{\rm th}$ root of unity. Let $A$ 
be a central simple algebra over $F$ of degree $n$. Suppose that $A$
is unramified on $R$ except possibly at  $x$ and $y$. Then   index$(A)
= $ index$(A \otimes F_{\nu})$   for the discrete valuation $\nu$ of $F$
given by the prime ideals either $(x)$ or $(y)$ of $R$. 

\paragraph*{Proof.}   
Suppose that $A$ is unramified  on $R$. Let $\nu$ be the discrete valuation of $F$ 
given by $x$.  Since $A$ is unramified 
on $R$,  by (2.2), we have  index$(A \otimes F_{\nu}) = $ index$(A)$.  
 
Suppose that $A$ is ramified on $R$ only at the prime ideal $(x)$.
Let $\nu$ be the discrete valuation on $F$ given by the prime ideal $(x)$ of $R$. 
By ([S1]), we have $A  = A' \otimes (u, x)$ for some unit $u$ in $R$ and
a central simple algebra $A'$ over $F$ which is unramified on $R$. 
Since $F_{\nu}$ is complete and $u$ is a unit at $\nu$, there is a unique extension of $\nu$ to 
$F_{\nu}(\sqrt[n]{u})$ such that $F_{\nu}(\sqrt[n]{u})$ is complete with the residue field 
$\kappa(\nu)(\sqrt[n]{\overline{u}})$, where $\overline{u}$ is the image of $u$ in $\kappa(\nu)$.
Since $A'$ is unramified on $R$,   we have 
$$
\begin{array}{rcll}
{\rm index}(A\otimes F_{\nu}) &  = & {\rm index}(A'\otimes (u, x) \otimes F_{\nu}) \\
& = & {\rm index}(A'\otimes F_{\nu}(\sqrt[n]{u})) \cdot [F_{\nu}(\sqrt[n]{u}) : F_{\nu} ] & ({\rm by ~(2.1)}) \\
& = & {\rm index}(A'\otimes F(\sqrt[n]{u})) \cdot [F_{\nu}(\sqrt[n]{u}) : F_{\nu} ] & ({\rm by ~(2.2)}) \\
& = & {\rm index}(A'\otimes F(\sqrt[n]{u})) \cdot [F(\sqrt[n]{u}) : F] & ({\rm by ~(2.3)}). 
\end{array}
$$
Since $A = A' \otimes (u, x)$,  the index of $A$ divides 
index$(A' \otimes F(\sqrt[n]{u})) \cdot [F(\sqrt[n]{u}) : F] = {\rm index}(A \otimes F_{\nu})$.  
Thus index$(A) = $ index$(A \otimes F_{\nu})$.  
  
Assume that $A$ is ramified on $R$  at both the primes $(x)$ and
$(y)$. Then by ([S1]), either $A  = A' \otimes (u_1, x) \otimes (u_2,
y)$ or $A  = A' \otimes (uy^r,  x)$ where $u_1, u_2, u$ are units in
$R$, $r$ coprime with $n$  and $A'$ unramified on $R$.

Suppose that $A  = A' \otimes (u_1, x) \otimes (u_2,  y)$ for some
units $u_1, u_2 \in R$ and $A'$ unramified on $R$. Let $\nu$
be the discrete valuation on $F$ given by $y$. Then by (2.1),
we have  index$(A \otimes F_{\nu}) = $ index$(A' \otimes (u_1, x) \otimes
F_{\nu}(\sqrt[n]{u_2})) \cdot [F_{\nu}(\sqrt[n]{u_2}) :
F_{\nu}]$. By (2.2), we have $[F_{\nu}(\sqrt[n]{u_2}) :
F_{\nu}] = [F(\sqrt[n]{u_2}) :
F]$. We now compute index$(A' \otimes (u_1, x) \otimes
F_{\nu}(\sqrt[n]{u_2}))$. 

Since $A' \otimes (u_1, x)$ is unramified at $\nu$ and
$u_2$ is a unit at $\nu$,   index$(A'\otimes(u_1, x) \otimes
F_{\nu}(\sqrt[n]{u_2})) =   {\rm index} (A'\otimes(u_1, x) \otimes
\kappa(\nu)(\sqrt[n]{\overline{u}_2}))$ (cf. 1.1), where $\overline{u_2}$ 
is the image of $u_2$ in $\kappa(\nu)$. 
Since $\kappa(\nu)$ is the  field of fractions of
$R/(y)$, the image of $x$ in $\kappa(\nu)$ gives a discrete valuation $\mu$
on the residue field $\kappa(\nu)$ with $\kappa(\mu) = k$ 
and  $\kappa(\nu)$ is complete with respect to $\mu$.  
Since $\overline{u_2}$ is a unit at $\mu$, 
$\kappa(\nu)(\sqrt[n]{\overline{u_2}})$ is a complete discrete valuated field 
with $\overline{x}$ as parameter and residue field
$k(\sqrt[n]{\overline{u_2}})$. 
Thus,  we have 
$$
\begin{array}{l}
{\rm index}(A'\otimes(u_1, x) \otimes
F_{\nu}(\sqrt[n]{u_2}))  =   {\rm index} (A'\otimes(u_1, x) \otimes
\kappa(\nu)(\sqrt[n]{\overline{u}_2}))  \hfill{({\rm by}~(1.1))} \\
  =  {\rm index}(A' \otimes
\kappa(\nu)(\sqrt[n]{\overline{u_2}},
\sqrt[n]{\overline{u_1}})) \cdot
[\kappa(\nu)(\sqrt[n]{\overline{u_2}},
\sqrt[n]{\overline{u_1}}) :
\kappa(\nu)(\sqrt[n]{\overline{u_2}})]   \hfill{({\rm by} ~(2.1))} \\
 =  {\rm index}(A' \otimes
F_{\nu}(\sqrt[n]{u_2},
\sqrt[n]{u_1})) \cdot
[F_{\nu}(\sqrt[n]{u_2},
\sqrt[n]{u_1}) :
F_{\nu}(\sqrt[n]{u_2})] \hfill{ ({\rm by}~(1.1) )}\\
 =  {\rm index}(A' \otimes
F(\sqrt[n]{u_2},
\sqrt[n]{u_1})) \cdot
[F(\sqrt[n]{u_2},
\sqrt[n]{u_1}) :
F(\sqrt[n]{u_2})] \hfill{({\rm by}~(2.2), (2.3)).}
\end{array}
$$
Hence 
$$
\begin{array}{l}
{\rm index}(A\otimes F_{\nu})  =  {\rm index}(A'\otimes(u_1, x) \otimes
F_{\nu}(\sqrt[n]{u_2})) \cdot [F_{\nu}(\sqrt[n]{u_2}) : F_{\nu}] \hfill{(by~(2.1))}\\
 =  {\rm index}(A' \otimes
F(\sqrt[n]{u_2}, \sqrt[n]{u_1})) \cdot
[F(\sqrt[n]{u_2}, \sqrt[n]{u_1}) : F(\sqrt[n]{u_2})] \cdot [F(\sqrt[n]{u_2}) : F] \\
 =  {\rm index}(A' \otimes
F(\sqrt[n]{u_2}, \sqrt[n]{u_1})) \cdot
[F(\sqrt[n]{u_2}, \sqrt[n]{u_1}) : F].
\end{array}
$$
On the other hand we have ${\rm index}(A) = {\rm index}(A'\otimes
(u_1, x) \otimes (u_2, y))$ divides ${\rm index}(A' \otimes
F(\sqrt[n]{u_1}, \sqrt[n]{u_2})) \cdot [F(\sqrt[n]{u_1},
\sqrt[n]{u_2}) : F]$.  In particular index$(A \otimes F_{\nu}) = {\rm index}(A)$.
 
Assume that $A = A' \otimes (uy^r,  x)$ for some unit $u
\in R$, $r$ coprime to $n$ and $A'$ unramified on $R$. Let
$\nu$ be the discrete valuation on $F$ given by the prime ideal $(x)$
of $R$.   
By  (2.1), we have index$(A \otimes F_{\nu}) = {\rm index}
(A'\otimes F_{\nu}(\sqrt[n]{uy^r})) \cdot
[F_{\nu}(\sqrt[n]{uy^r}) : F_{\nu}].$  Since $y$ is a regular prime in $R$ and 
$r$ is coprime to $n$, it follows that $[F(\sqrt[n]{uy^r} ): F] = 
[F_{\nu}(\sqrt[n]{uy^r} ): F_{\nu}] = n$.   
Since $F_{\nu}$ is a complete discrete valuated field and $A'$ is unramified,
we have ${\rm index}
(A'\otimes F_{\nu}(\sqrt[n]{uy^r})) =  {\rm index}(A' \otimes
\kappa(\nu)(\sqrt[n]{\overline{uy^r}}))$ (cf. 1.1).
Since $\kappa(\nu)$ is a complete discrete valuated field with $\overline{y}$ as a parameter and residue field $k$, 
$\kappa(\nu)(\sqrt[n]{\overline{uy^r}})$ is also a complete discrete valuated field with residue field $k$. 
Since $A'$ is unramified on $R$, we have  ${\rm index}(A' \otimes \kappa(\nu)(\sqrt[n]{\overline{uy^r}})) = 
{\rm index}(A'\otimes k)$ (cf. (1.1)). Since $R$ is complete, by
(1.1),  we have $ {\rm index}(A') = {\rm index}(A'\otimes k) =  {\rm index}(A' \otimes
\kappa(\nu)(\sqrt[n]{\overline{uy^r}})) = {\rm index}(A' \otimes F_{\nu}(\sqrt[n]{uy^r})).$ 
Therefore, we have
$$
\begin{array}{rcl}
{\rm index}(A \otimes F_{\nu}) & = & {\rm index}(A' \otimes (uy^r, x) \otimes F_{\nu}) \\
& = & {\rm index}(A ' \otimes F_{\nu}(\sqrt[n]{uy^r})) \cdot [F_{\nu}(\sqrt[n]{uy^r}) : F_{\nu}] \\
& = & {\rm index}(A' ) \cdot n.
\end{array}
$$
Since $A = A' \otimes (uy^r, x)$, index$(A) \leq {\rm index}(A') \cdot
n = {\rm index}(A \otimes F_{\nu})$. 
Thus index$(A) = {\rm index}(A \otimes F_{\nu})$.   \hfill $\Box$

\paragraph*{Proposition 2.5.} Let $R$ be a 2-dimensional regular ring with field of
fraction $F$ and $n$ an integer which is a unit in $R$. Assume that $F$ contains a primitive 
$n^{\rm th}$ root of unity.  
Suppose that $R$ is complete with respect to $(s)$-adic topology for some
prime element  $s \in R$ with $R/(s)$ a Dedekind domain. Let $\nu$ be the discrete valuation 
on $F$ given by $s$. Let $A$ be a central simple 
algebra over $F$ of degree $n$ which is unramified on $R$ except at $\nu$. Further assume that
the residue of $A$ at $(s)$ is given by a unit $a$ in $R/(s)$. 
Then index$(A) =$ index$(A \otimes F_{\nu})$. 

\paragraph*{Proof.} Suppose that $A$ is unramified on $R$.  
Since $R$ is $(s)$-adically complete and $R/(s)$ is a regular domain
of dimension 1  with field of fractions $\kappa(\nu)$,  index$(A) = $
index$(A \otimes _{R/(s)} \kappa(\nu)) = $ index$( A \otimes F_{\nu})$
(cf. (1.1)). 

Suppose that $A$ is ramified on $R$. Then by the assumption on  $A$, $A$ is
ramified on $R$ only at the prime ideal $(s)$ of $R$ and the residue of $A$
at $(s)$ is given by a unit $a$ in $R/(s)$.  
Let $u \in R$ with image $a \in R/(s)$. Since $R$ is
$(s)$-adically complete and the image $a$ of $u$ modulo $(s)$ is a
unit, $u$ is a unit in $R$.  
The cyclic algebra $(u, s)$
is ramified on $R$ only at $(s)$ and the residue of $(u,
s)$ at $(s)$ is $(a)$.  Since $A$ is ramified only at the prime
ideal $(s)$ of $R$, and $(a)$ is the residue of $A$ at $s$, we
have $A  = A' \otimes (u, s)$ for some central simple algebra
$A'$ over $F$ which is unramified on $R$. By (2.1), we have index$(A \otimes F_{\nu}) = 
{\rm index}(A' \otimes F_{\nu}(\sqrt[n]{u})) \cdot [F_{\nu}(\sqrt[n]{u}) : F_{\nu}]$. Since 
$R$ is $(s)$-adically complete and $R/(s)$ is integrally closed
domain with  $u$ a unit, 
we have $[F(\sqrt[n]{u}) : F] = [F_{\nu}(\sqrt[n]{u}) : F_{\nu}]$.
Since $A'$ is unramified on $R$  and $R$ is
$(s)$-adically complete, by (2.2), we have 
 index$(A'\otimes F_{\nu}(\sqrt[n]{u})) = {\rm index}(A' \otimes F(\sqrt[n]{u}))$.
Hence we have index$(A \otimes F_{\nu}) = {\rm index}(A' \otimes F(\sqrt[n]{u})) \cdot [F(\sqrt[n]{u}) : F] $.
In particular index$(A \otimes F_{\nu}) = {\rm index}(A)$. \hfill $\Box$

\paragraph*{Theorem 2.6.} Let $T$ be
a complete discrete valuation ring with fraction field $K$ and
residue field $k$. Let $X$ be a regular, projective, geometrically
integral curve over $K$ and $F=K(X)$  the function field of $X$. Let
$l$ be prime not equal to  char$(k)$ and $A$  a central simple algebra
over $F$ of index  $n = l^d$ for some $d \geq 1$. 
Assume that $K$ contains a primitive $n^{th}$
root of unity. Then index$(A) = $ index$(A \otimes F_{\nu})$  for some discrete
valuation ${\nu}$ of $F$. 

\paragraph*{Proof.} Let $A$ be a central simple algebra over $F$.
We choose a regular proper model ${\cal X}/T$ of $X/K$  such that
the support of the ramification divisor $A$ and the components of
the special fibre of ${\cal X}/T$ are a union of regular curves with
normal crossings. Let $Y={\cal X} \times_{T}k$ denote the special
fibre.

For each irreducible curve $C$ in the support of ramification
divisor of $A$, let $a_C \in \kappa(C)^*$ be such that the residue
of $A$ at $C$ is $(a_C) \in H^1(\kappa(C)^*, \mathbb{Z}/n\mathbb{Z})$.

Let $S$ be a finite set of closed points of the special fibre
containing all singular points of $Y$ and all those points of irreducible
curves $C$ where $a_C$ is not a unit.
 
For each $P \in S$,  let $F_P$ be the field of fractions of the
completion $\hat{R}_P$ of the local ring $R_P$ of ${\cal X}$ at $P$.
Let $t$ be a uniformizing parameter for $T$.
For each  irreducible component $U$ of $ Y\setminus S$, let $R_U$ be the ring of
elements in $F$ which are regular on $U$. It is a regular ring. 
Let $\hat{R}_U$ be the completion of $R_U$ at $(t)$ and $F_U$ the 
field of fractions of  $\hat{R}_U$. As in ([CTPS], proof of 3.1), we
choose $S$ such that  for every irreducible component $U$ of $Y
\setminus S$,   $t=u.s^r$ for some integer $r \geq 1$, a unit $u
\in R_{U}$ and  $R_{U}/s= {\hat R}_U/s$  is a
Dedekind domain with field of fractions $k(U)$.  In particular, 
the $t$-adic completion ${\hat R}_U$ coincides with
the $s$-adic completion of $R_{U}$.

Let ${\cal U}$ be the set of irreducible components of $X \setminus S$.
By ([HHK1], 5.1),  we
have index$(A) =$ lcm of index$(A \otimes F_{\zeta})$, $\zeta \in {\cal U}
\cup S$. Since index$(A \otimes F_{\zeta})$ is a power of the prime
$l$ for all $\zeta \in {\cal P} \cup P$, we have  index$(A) =$
index$(A \otimes F_{\zeta})$  for some $\zeta \in {\cal U }
\cup S$. In particular index$(A)$ equal to either index$(A \otimes
F_U)$  for some irreducible components $U$ of
$Y \setminus S$ or index$(A \otimes F_P)$   for some $P \in S$.

Suppose that index$(A) = $ index$(A \otimes F_P)$   for some $P \in S$. The
local ring $\hat{R}_P$ is a complete regular local ring of dimension 2 with maximal ideal  $(x,
y)$ such that $A $ is ramified on $\hat{R}_P$ at most at $x$ and $y$. By
(2.4),  we have index$(A) = $ index$(A \otimes {F_P}_{\nu})$   for the discrete valuation
of $F_{P}$ given by either $(x)$ or $(y)$.  Since the  restriction
$\nu_0$  of $\nu$ to $F$ is non-trivial,
  $\nu_0$   is a discrete valuation on $F$ and index$(A) = $ index$(A \otimes F_{\nu_0})$. 
 
Suppose that index$(A) = $ index$(A \otimes F_U)$   for some irreducible
component $U$ of $X \setminus S$. Then, by the choice of $U$, 
$A$ is unramified on $\hat{R}_U$ except at $(s)$ and the residue at $(s)$ is  
given by a unit in $R_U/(s)$. Let $\nu$ be the discrete valuation on $F_U$ given by $(s)$.
Since $\hat{R}_U$ is $(s)$-adically complete,  by (2.5),  index$(A) =
$ index$(A \otimes {F_U}_{\nu})$.
Since the restriction of $\nu$ to $F$ is also given by the ideal $(s)$ in $R_U$, $\nu$ is non-trivial on $F$.
In particular $F_{\nu} \subset {F_U}_{\nu}$ and index$(A) = $ index$(A \otimes F_{\nu})$.  \hfill $\Box$.

\paragraph*{Remark 2.7.}  Let $F$ and
$A$ be as in the above theorem. Then by (2.6), it follows that there
exists a discrete valuation $\nu$ of $F$ such that index$(A) =$
index$(A \otimes F_{\nu})$.  From the proof of (2.6) it follows that this
discrete valuation comes from a codimension one point of a regular
proper model of $F$. In particular, the residue field $\kappa(\nu)$
is either a finite extension of $K$ or a function field of a curve
over a finite extension of $k$.

\section*{3. Necessary conditions for Admissibility }

In this section, we give a necessary condition for a finite group to be
admissible over function fields of curves over complete discretely
valued fields.

Let $K$ be a complete discretely valued field with residue field $k$ and $\pi$ a parameter.
Let $D$ be a central division algebra over $K$ of degree $n$ and  $L$  a maximal subfield of $D$. 
Suppose that   $n$ is coprime to char$(k)$. Let $(E_0, \sigma_0) \in H^1(k,
\mathbb{Z}/n\mathbb{Z})$ 
be the  residue of $D$ and $(E, \sigma) \in H^1(K, \mathbb{Z}/n\mathbb{Z})$ 
the lift of $(E_0, \sigma_0)$ (cf. \S 1).  We fix an algebraic closure $\overline{K}$ of $K$ and 
also an algebraic closure $\overline{k}$ of $k$.
All the finite extensions of $K$ and $k$  are considered as subfields of $\overline{K}$ and 
$\overline{k}$ respectively. 
Let $F = L \cap E$.  Let $L_1$ be the maximal unramified extension 
of $F$ contained in $L$. Since $E/K$ is unramified, $L \cap E = L_1 \cap E = F$.

\paragraph*{Lemma 3.1.}   Let $K$, $D$, $L$, $L_1$, $(E_0, \sigma_0)$, $(E, \sigma)$, and $F$
be as above. Then $D \otimes_K F \otimes_F (E/F, \sigma^{[F:K]}, \pi)^{op}$
is unramified at the discrete valuation of $F$.

\paragraph*{Proof.}  Since the residue of $D$ is $(E_0, \sigma_0)$, $D
\otimes (E, \sigma, \pi)^{op}$ is unramified at the discrete valuation
of $K$.  In particular, $D
\otimes _K (E, \sigma, \pi) \otimes F$ is unramified at the discrete
valuation of $F$.  Since $F \subset
E$,  we have $(E,\sigma, \pi)  \otimes_K F = (E/F, \sigma^{[F :
  K]}, \pi)$ in $Br(F)$. Thus  $(D \otimes_K F) \otimes_F (E/F, \sigma^{[F:K]}, \pi)^{op}$
is unramified at the discrete valuation of $F$. \hfill $\Box$

\paragraph*{Lemma 3.2.}  Let $K$, $D$, $L$, $L_1$,  $(E_0, \sigma_0)$, $(E, \sigma)$ and $F$
be as above. Then $[E : F]$ divides $[L : L_1]$. 

\paragraph*{Proof.}   We have the following commutative diagram (cf.,  [S1],
10.4)
\begin{displaymath}
\xymatrix{ _nBr(K) \ar[r] \ar[d]_{res} &
H^{1}(k,\mathbb{Z}/n\mathbb{Z}) \ar[d]^{e.res} \\ _nBr(L)  \ar[r]  &
H^{1}(L_0,\mathbb{Z}/n\mathbb{Z}) },
\end{displaymath}
where $L_0$ is the residue field of $L$ and  $e$ the ramification index of $L/K$.
From the above commutative diagram, the residue of
$D \otimes L$ is the restriction of $e(E_0, \sigma_0)$ to $L_0$.  
Since $L$ is a maximal subfield of $D$, $D \otimes L$ is a split algebra.
In particular, the residue of $D \otimes L$ is trivial. Hence the restriction of $e(E_0, \sigma_0)$ to $L_0$ is
trivial and $[E : F]$ divides $e $  (cf. \S 1).   Since $[E : F] =
[E_0 : F_0]$ and $e = [L : L_1]$, $[E : F]$ divides $[L : L_1]$.  \hfill $\Box$

\paragraph*{Lemma 3.3.}  Let $K$, $D$, $L$, $L_1$, $(E_0, \sigma_0)$, $(E, \sigma)$ and $F$
be as above.  Then  index$(D \otimes_K E) = [L_1E : E]$.

\paragraph*{Proof.} Since $(E/F, \sigma^{[F:K]}, \pi)^{op} \otimes E$ is split,
we have $D \otimes_K E = D \otimes_KF \otimes_F (E/F, \sigma^{[F:K]}, \pi)^{op} \otimes_F E$.
By (3.1), $D \otimes_K F \otimes_F (E/F, \sigma^{[F:K]}, \pi)^{op}$ is unramified at the discrete 
valuation of $F$.  Hence $D \otimes_K E$ is unramified at the discrete valuation of $E$. 
Since $L$ is a maximal subfield of $D$, $D \otimes_K L$ is
split. Since $L/L_1$ is totally ramified, $LE/L_1E$ is totally
ramified. In particular, the residue field of $LE$ and $L_1E$ are
equal.  Since  $D \otimes_K E$ is unramified and $(D \otimes_KE)
\otimes_E LE$ is plit,  $(D \otimes_K E) \otimes_E L_1E$ is split
(cf. 1.1).  In particular  index$(D \otimes_K E) \leq
[L_1E : E]$.  

Since   $(D \otimes_K F) \otimes _F (E/F, \sigma^{[F:K]}, \pi)^{op}$ is unramified and $\pi$ is a
parameter in $E$, by (2.1),  we have 
$$ 
\begin{array}{l}
{\rm index}(D \otimes_K F) = {\rm index}( (D \otimes_K F) \otimes _F (E/F, \sigma^{[F:K]}, \pi)^{op} \otimes_F 
(E/F, \sigma^{[F:K]}, \pi))  \\
= {\rm index}((D \otimes_K F) \otimes_F (E/F, \sigma^{[F:K]}, \pi)^{op}
\otimes_F E ) \cdot [E
: F] \\
= {\rm index}(D\otimes_K E) \cdot [E : F] 
\end{array}
$$
On the other hand, since $F \subset L$ and $L$ is a maximal subfield
of $D$,  index$(D \otimes_K F) = [L : F]$.
Hence   ${\rm index}(D \otimes_K E) \cdot  [E : F] = [L : F] = [L :
L_1]\cdot [L_1 : F]$.
By (3.2), we have  $[E : F]$ divides $  [L : L_1]$. Hence 
${\rm index}(D \otimes_K E)  \geq  [L_1 : F] \geq [L_1E : E]$. Therefore
index$(D \otimes_K E) = [L_1E : E]$.
\hfill $\Box$

\paragraph*{Lemma 3.4.} Let $K$ be a complete discretely valued
field with residue field $k$ and $P$  a $p$-group, $p$ a
prime. Suppose that  $p$ is coprime to   char($k)$. 
If $P$ is admissible over $K$, then $P$ has a
normal series   $P\supseteq P_{1}\supseteq P_{2}$ such that
 
(1) $P/P_{1}$ and $P_{2}$ are cyclic

(2) $P_{1}/P_{2}$ is admissible over some finite extension of  $k$.

\paragraph*{Proof} Suppose that $P$ is admissible over $K$.
Then there exists a Galois extension $L/K$ and a division ring $D$
central over $K$ which contains $L$ as maximal subfield such that
$P=G(L/K)$. Let $L_0$ be the residue field of $L$.   Let
$\partial(D)$=$(E_{0},\sigma_{0})\in
H^{1}(k,\mathbb{Z}/n\mathbb{Z})$ be the residue of D and
$(E,\sigma)\in H^{1}(K,\mathbb{Z}/n\mathbb{Z})$ be the lift of
$(E_{0},\sigma_{0})$.    
Let $L_1$ be the maximal unramified extension of $K$
contained in $L$. 
Since $E$ is unramified extension of $K$,
we have $L \cap E = L_1 \cap E$.  Let $F = L\cap E$. Since $E/K$
is cyclic, $F/K$ is also cyclic.

Let $P_1$ be the Galois group of $L/F$.  Since $F/K$ is cyclic,
$P_1$ is a normal subgroup of $P$ and $P/P_1$ is cyclic.

Let $P_2$ be the Galois group of $L/L_1$. Then $P_2$ is a
subgroup of $P_1$. Since $L/L_1$ is a totally ramified Galois
extension of degree coprime to char$(k)$,  $P_2$ is cyclic ([Se],
Cor.2  and Cor.4 of Ch.IV, \S 2). Since $L_1/F$ is a Galois
extension (cf. \S 1), $P_2$ is a normal subgroup of $P_1$. The residue field
$F_0$ of $F$ is the same as the intersection of $L_0$ and $E_0$. We now show
that $P_1/P_2$ is admissible over  $E_0$.

Let $D'$ be the division algebra with center $E$ which is Brauer
equivalent to  $D\otimes _KE $.  
Since $D' \otimes_EL_1E =  (D \otimes_K E) \otimes_E L_1E $ and $(D
\otimes_K E) \otimes_E L_1E$ is split
(cf. proof of (3.3)), $D' \otimes_E L_1E$ is also split. Since, by (3.3)  degree$(D') = [L_1E : E]$, 
$L_1E$ is a maximal subfield of $D'$.  By (3.1), $D' = D \otimes_K E$
is unramified at the discrete valuation of $E$.  Let  $\overline{D'}$ be the  image of $D'$
over the residue field $E_0$ of $E$. Since $E$ is complete, $\overline{D'}$ is central
division algebra over $E_0$ and $L_0E_0$ is a maximal subfield of
$\overline{D'}$. Hence   $Gal(L_0E_0/E_0)$ is admissible over $E_0$.
 Since the residue field of $L_1E$ is $L_0E_0$ and
$L_1E/E$ is an unramified Galois extension, we have $Gal(L_0E_0/E_0)
\simeq Gal(L_1E/E)$.
Since $L_1/K$ and $E/K$ are Galois and $F = L_1 \cap E$,  we have $Gal(L_1E/E) \simeq
Gal(L_1/F)$.  Since  $ P_1/P_2 \simeq Gal(L_1/F) \simeq
Gal(L_0E_0/E_0)$, $P_1/P_2$ is admissible over $E_0$.
\hfill $\Box$ 

\vskip 3mm

The above lemma  immediately  gives the following

\paragraph*{Proposition 3.5.} Let $K$ be a complete discretely valued
field with residue filed $k$ and $G$ be a finite group. Suppose that
 the order of $G$ is coprime to char$(k)$.
If $G$ is admissible over $K$ then
every Sylow subgroup $P$ of $G$ has a normal sereis   $P\supseteq
P_{1}\supseteq P_{2}$ such that

(1) $P/P_{1}$ and $P_{2}$ are cyclic

(2) $P_{1}/P_{2}$ is admissible over some finite extension of the
residue field of $K$.

\paragraph*{Proof.} Let $G$ be an admissible group over $K$.
Then there is  a field extension $L/K$ and a division algebra $D$
central over $K$ containing $L$ as a maximal subfield with
$Gal(L/K)=G$. Let $P$ be a Sylow subgroup of $G$. Let $L^{P}$
be the fixed of $P$. Then $L^{P}$ is a complete discretely valued
field. Let $D'$ be the commutant of $L^{P}$ in $D$. Then $D'$ is a
central division algebra over $L^{P}$ and $G(L/L^{P})=P$ is
admissible over $L^{P}$. Since $L^{P}$ is also a complete discrete
valued field, the result follows by (3.4). \hfill $\Box$

\paragraph*{Corollary 3.6.} Let $K$ be a complete discretely valued field
with residue filed $k$ either a local field or a global field. Let
$G$ be a finite group such that char($k$) coprime to the order of
$G$.  If $G$ is admissible over $K$ then every Sylow subgroup $P$ of $G$ has
a normal series  $P\supseteq P_{1}\supseteq P_{2}$ such
that
 
(1) $P/P_{1}$ and $P_{2}$ are cyclic

(2) $P_{1}/P_{2}$ is  metacyclic.

\paragraph*{Proof.} Every finite extension of the residue field is
either a local field or a global field. The corollary follows from
(3.5) and ([Sc]). \hfill $\Box$

\paragraph*{Theorem 3.7.} Let $K$ be a complete discretely valued
field with residue field $k$. Let  $F$ be the function field of a curve
over  $K$.  Let $G$ be a
finite group  with order coprime to char$(k)$. If $G$ is admissible over
$F$ then every Sylow subgroup $P$ of $G$ has a normal series  
$P\supseteq P_{1}\supseteq P_{2}$ such that
 
(1) $P/P_{1}$ and $P_{2}$ are cyclic

(2) $P_{1}/P_{2}$ is admissible over some finite extension of the
residue field at a discrete valuation of $F$.

\paragraph*{Proof}  First we reduce to a Sylow subgroup as in (3.5).
Let $G$ be an admissible group over $F$.
Then there is  a field extension $L/F$ and a division algebra $D$
central over $F$ containing $L$ as a maximal subfield with
$Gal(L/F)=G$. Let $P$ be a Sylow subgroup of $G$. Let $L^{P}$
be the fixed of $P$. Then $L^{P}$ is a complete discretely valued
field. Let $D'$ be the commutant of $L^{P}$ in $D$. Then $D'$ is a
central division algebra over $L^{P}$ and $G(L/L^{P})=P$ is
admissible over $L^{P}$. Since $L^P$ is a finte extension of $F$,
$L^{P}$ is also a function field of a curve over a finite extension of
$K$. Since the degree of $D'$ is a power of prime and coprime to the
char$(k)$, by (2.6),  there exists a discrete valuation $\nu$ of $L^P$
such that $D' \otimes_
{L^P} L^P_{\nu}$ is   division. Since $L\otimes_{L^P} L^P_{\nu}$ is a maximal
subfield of $D' \otimes_{L^P} L^P_{\nu}$ and $P = Gal(L \otimes_{L^P}
L^P_{\nu}/L^P_{\nu})$,  the result follows from  (3.5).  \hfill $\Box$

\vskip 3mm

The following is immediate from (3.7) and (3.6).

\paragraph*{Corollary 3.8.} Let $K$ be a local adic field and F the function
field of a curve over  $K$.
Let $n$ be a natural number which is coprime to the characteristic of
the residue field of $K$  and $G$ a finite
group of order $n$. If $G$ is admissible over $F$ then every Sylow
subgroup $P$ of $G$ has a normal series  $P\supseteq P_{1}\supseteq
P_{2}$ such that
 
(1) $P/P_{1}$ and $P_{2}$ are cyclic

(2) $P_{1}/P_{2}$ metacyclic.

\vskip 3mm

The following is proved in ([HHK2], 4.5).

\paragraph*{Corollary 3.9.} Let $K$ be a complete discretely valued
field with residue field algebraically closed. Let F be the function
field of a curve over $K$. Let $G$ be a finite group of order $n$
such that the characteristic of the residue field of $K$ is coprime
to $n$. If $G$ is admissible over $F$ then every Sylow subgroup
$P$ of $G$  is metacyclic.

\paragraph*{Proof.} By (3.7), every Sylow subgroup $P$ of $G$ has a
normal series  $P\supseteq P_{1}\supseteq P_{2}$ such that
 
(1) $P/P_{1}$ and $P_{2}$ are cyclic

(2) $P_{1}/P_{2}$ is admissible over some finite extension of the
residue field of $K$.

Let $k$ be the residue field of $K$.  By (2.7) and the proof of
(3.7), the residue field at the discrete valuation given in (2) is
either a finite extension $K$ or the function field of a curve over
$k$. Since $k$ is algebraically closed, these residue fields have
cohomological dimension at most one and there are no non-trivial
division algebras over such fields. Hence $P_1 = P_2$. \hfill $\Box$

\vskip 3mm

We now give an example of a finite group which is not
$\mathbb{Q}_{p}(t)$-admissible.

\paragraph*{Example 3.10.} Let $\ell$ and $p$ be two distinct primes.
Let $P =(\mathbb{Z}/l\mathbb{Z})^5$. We claim that this group is not
admissible over $\mathbb{Q}_{p}(t)$. Suppose that $P$ is admissible
over $\mathbb{Q}_{p}(t)$. Then by (3.8), there is a normal series   $P\supseteq
P_{1}\supseteq P_{2}$ such that
 
(1)$P/P_{1}$ and $P_{2}$ are cyclic

(2)$P_{1}/P_{2}$ is metacyclic

Since $P_{2}$ and $P/P_{1}$ are cyclic their orders will be at most
$l$.  This implies that $|P_{1}/P_{2}|\geq l^3$.  Since every
element of $P_{1}/P_{2}$ has order at most $l$  and
$|P_{1}/P_{2}|\geq l^3$, $P_{1}/P_{2}$ cannot be metacyclic.
\hfill $\Box$

\paragraph*{Remark 3.11.} Let $F$ be a field and $p$ a prime.  Suppose
that for every  finite extension $E$ of $F$
there is a set $\Omega_E$ of discrete
valuations of $E$ such that given a central division algebra $D$ over
$E$ of degree a power of $p$, there exists a discrete valuation $\nu
\in \Omega_E$  such that $D
\otimes E_{\nu}$ is division. Let $G$ be a finite group.  
The proof of (3.7) gives us the following: 
If $G$ is admissible over $F$, then every  
$p$-Sylow subgroup of $G$ has filtration as in (3.7).

\paragraph*{Remark 3.12.}    Let $G$ be a finite group satisfying the
conditions of (3.5). Then   every homomorphic image of
$G$ also satisfy the same conditions. However there  is an  example of
a group $G$ with a
homomorphic image $H$ such that $G$ is admissible over the complet
discrete valuation field $\mathbb{Q}((t))$ but $H$ is not admissible
([FS]).  Hence the conditions given in (3.5) for a group to be admissible
are not sufficient.

\section*{4. A class of Admissible groups over
$\mathbb{Q}_{p}(t)$}

Let $K$ be a discretely valued field with residue field $k$. Let $F$
be the function field of a curve over $K$. Let $n$ be an integer
which is coprime to the characteristic of $k$. Suppose that $K$
contains a primitive $n^{th}$ root of unity. Then in ([HHK2], 4.4) it is
proved that every finite group of order $n$ with every Sylow
subgroup product of at most two cyclic groups  is admissible over $F$. They used
the patching techniques to prove this result. In this section we
prove a similar result for groups with every Sylow subgroup is a product of
at most 4 cyclic  groups,  with an additional assumption on the residue field
$k$. We begin with the following

\paragraph*{Lemma 4.1.} Let $R$ be a regular local ring of
dimension two with residue field $k$ and field of fraction $F$. Let
$n_1$ and $n_2$ be  natural numbers  which are coprime to the
char$(k)$.  Assume that $F$ contains a primitive $(n_1n_2)^{th}$ root
of unity and there is an element in $k^*/k^{*n_2}$ of order $n_2$.
Then there is a central division algebra $D$ over $F$ of degree
$n_1n_2$.

\paragraph*{Proof.} Let $m$ be the maximal ideal of $R$.  Since $R$ is
a regular local ring of dimension two, 
we have $m = (t, s)$. By the assumption on $k$, there is an element
$\lambda_0 \in k^*$ such that its order in $k^*/k^{*n_2}$ is $n_2$.
Let $\lambda \in R$ which maps to $\lambda_0$. Let $a \in R$ be a
unit with $a^{n_1} \neq 1$. Let $\xi_1$ be  a primitive $n_1^{th}$
root of unity and $\xi_2$   a primitive $n_2^{th}$ root of unity.

Let
$$D_1 = (\frac{s}{s-t}, \frac{s-t^2}{s-a^{n_1}t^2})_{n_1}$$
and
$$
D_2 = (\frac{s}{s-t^2}, \frac{s-\lambda t^2}{s-t^2})_{n_2}.
$$
Let $D = D_1 \otimes_F D_2$. Then the degree of $D$ is $n_1n_2$.
We now show that $D$  is a division algebra.

Let $S = R[x]/(s-t^2x)$. Then the field of fractions of $S$ is
isomorphic to $F$. We have
$$
D_1 = (\frac{tx}{tx-1}, \frac{x-1}{x-a^{n_1}})_{n_1}
$$
and
$$
D_2 = (\frac{x}{x-1}, \frac{x-\lambda}{x-1})_{n_2}.
$$

The ideal $(t)$ of $S$ is a prime ideal and gives a discrete
valuation $\nu$ on $F$. Let $ F_{\nu}$ be the completion of $F$ at
$\nu$. To show that $D_1 \otimes_F D_2$ is a division algebra, it is
enough to show that $D_1 \otimes_F D_2 \otimes _FF_{\nu}$ is a division
algebra. Since $tx/tx-1$ is a parameter at $\nu$, by (2.1), we have
$$
{\rm index}(D_1 \otimes_F D_2 \otimes_F F_{\nu}) = {\rm index}(D_2 \otimes_F
F_{\nu}(\sqrt[n_1]{\frac{x-1}{x-a^{n_1}}})) \cdot
[F_{\nu}(\sqrt[n_1]{\frac{x-1}{x-a^{n_1}}}) : F_{\nu}].
$$
Since $\frac{x-1}{x-a^{n_1}}$ is a unit at $\nu$ and the residue
field $\kappa(\nu)$ at $\nu$ is $k(x)$, by the assumption on $a$, we have
$[F_{\nu}(\sqrt[n_1]{\frac{x-1}{x-a^{n_1}}}) : F_{\nu}] = n_1$.
Since $D_2$ is unramified at $\nu$, the index of $D_2 \otimes_F
F_{\nu}(\sqrt[n_1]{\frac{x-1}{x-a^{n_1}}})$ is equal to the index of
its image $(\frac{x}{x-1}, \frac{x-\lambda}{x-1})_{n_2}$ over the
residue field $k(x)(\sqrt[n_1]{\frac{x-1}{x-a^{n_1}}})$. Let $\theta$
be the discrete valuation on $k(x)$ given by $(x)$. Let $\upsilon$ be
the extension of $\theta$ to $k(x)
(\sqrt[n_1]{\frac{x-1}{x-a^{n_1}}})$.  Then the residue
field of $\upsilon$ is $k$.  The residue of $(\frac{x}{x-1}, \frac{x-\lambda}{x-1})_{n_2}$ at $\upsilon$ is
the class of $\lambda_0$. Since the    order of
$\lambda_0$ in $k^*/k^{*n_2}$ is $n_2$, the index of $D_2$ is $n_2$.
Hence the index of $D_1 \otimes_F D_2 \otimes F_{\nu}$ is $n_1n_2$.
Since the degree of $D_1 \otimes_F D_2$ is $n_1n_2$, $D_1 \otimes_F D_2$
is a division algebra. \hfill $\Box$

\paragraph*{Theorem 4.2.} Let $K$ be a discretely valued field with residue field
$k$ and   $F$  the function field of a curve over $K$. Let $n$ be
an integer which is coprime to the characteristic of $k$. Suppose
that $K$ contains a primitive $n^{th}$ root of unity. Assume that
for every finite extension $L$ of $k$, there is an element in
$L^*/L^{*n}$ of order $n$. If $G$ is a finite group of order $n$ with every
Sylow subgroup is a quotient of $\mathbb{Z}^4$, then $G$ is admissible
over $F$.

\paragraph*{Proof.} Let $R$ be the ring of integers in $K$. Let
${\cal X}$ be a regular proper two dimensional scheme over $R$ with
function field $F$ and the reduced special fibre is a union of
regular curves with normal crossings. Let $p_1, \cdots, p_r$ be the
prime factors of $n$. Let $Q_1, \cdots, Q_r$ be regular closed
points on the special fibre of ${\cal X}$. Let $R_{Q_i}$ be the
regular local ring at $Q_i$, $\hat{R}_{Q_i}$ be the completion of
$R_{Q_i}$ at the maximal ideal and $F_{Q_i}$ the field of fractions
of $\hat{R}_{Q_i}$. Let $t_i \in R_{Q_i}$ be a prime defining the
irreducible component of the special fibre of ${\cal X}$ containing
$Q_i$. Let $P_i$ be a $p_i$-Sylow subgroup of $G$. By ([HHK2], 4.2), it is
enough to show that there exists a central division algebra $D_i$
over $F_{Q_i}$ and maximal subfield $L_i$ of $D_i$ with
$Gal(L_i/F_{Q_i}) \simeq P_i$ and $L_i \otimes \hat{F}_{Q_i}$ a
split algebra.

For a given $i$, $1 \leq i \leq r$, let $Q = Q_i$, $t = t_i$ and $P =
P_i$.  Since the residue field $\kappa(Q)$ of $\hat{R}_{Q}$ is a finite
extension of $k$, by the assumption on $k$, there is an element in
$\kappa(Q)^*/\kappa(Q)^{*n}$ of order $m$ for any $m$ dividing $n$. 
Since $P$ is a quotient of $\mathbb{Z}^4$,  $ P \simeq C_{n_1} \times C_{n_2} \times
C_{n_3} \times C_{n_4}$. Since
$\hat{R}_{Q}$ is regular local ring of dimension 2 and $t$ is a
regular prime, we have $m_{Q} = (t, s)$. Let $\xi_1$ be a
primitive ${n_1n_2}^{th}$ root of unity and $\xi_2$  a primitive
${n_3n_4}^{th}$ root of unity.

Let
$$D_1 = (\frac{s}{s - t }, \frac{s - t^2}{s - a^{n_1n_1}t^2})_{n_1n_2}$$
and
$$
D_2 = (\frac{s}{s - t^2}, \frac{s-\lambda
t^2}{s - t^2})_{n_3n_4}.
$$
for suitable $a$ and $\lambda$ as in (4.1). Let $D = D_1 \otimes_F
D_2$. Then, by (4.1), $D$ is a division algebra over $F_Q$.

In particular  $D_1$  and $D_2$ are division algebras over $F_{Q}$. The cyclic
algebra $D_1$ is generated by $x_1$ and $y_1$ with relations
$$x_1^{n_1n_2} = \frac{s}{s-t}, ~ y_1^{n_1n_2} =
\frac{s-t^2}{s-a^{n_1n_2}t^2} ~{\rm and} ~x_1y_1 =
\xi_1y_1x_1.
$$
Similarly $D_2$ is  generated by $x_2$ and $y_2$ with relations
$$x_2^{n_3n_4} = \frac{s}{s-t^2}, ~ y_2^{n_3n_4} =
\frac{s-\lambda t^2}{s-t^2} ~{\rm and} ~x_2y_2 =
\xi_2y_2x_2.
$$
Let $L_1$ be the subalgebra of $D_1$ generated by $x_1^{n_1}$ and
$y_1^{n_2}$.  Let $L_2$ be the subalgebra of $D_2$ generated by
$x_2^{n_3}$ and $y_2^{n_4}$. Then $ L = L_1 \otimes L_2$ is a
maximal subfield of $D_1 \otimes D_2$,  $Gal(L/F) = C_{n_1} \times
C_{n_2} \times C_{n_3} \times C_{n_4}$ and $L \otimes \hat{F}$ is a
split algebra (cf., ([HHK2]) proof of 4.4, ).  \hfill $\Box$

\paragraph*{Corollary 4.3.} Let $K$ be a local field  and $k$ its residue field.
Let  $F$ be the function field of a curve over $K$. Let $n$ be an
integer which is coprime to the characteristic of $k$. Suppose that
$K$ contains a primitive $n^{th}$ root of unity.   If $G$ is a group
of order $n$ with every Sylow subgroup  a quotient of $\mathbb{Z}^4$,
then $G$ is admissible over $F$.

\paragraph*{Proof.} Since $k$ is a finite field, for any finite
extension   $L$  of $k$ and for any
natural number $n$ coprime to the characteristic of $L$, we have
$L^*/L^{*n}$ is cyclic group of order $n$. Hence $k$ satisfies the
condition of (4.2). \hfill $\Box$.

\vskip 3mm
 
Let $K$ be a complete discretely valued field with residue field
$k$. Let $G$ be a finite abelian group of order $n$ which is a product
of at most four cyclic groups. Suppose that  $n$ is coprime to char$(k)$ and the
order of $k^*/k^{*n}$ is at least $n$. 
If $K$ contains a primitive $n^{th}$ root of unity, then by (4.3), $G$ is
admissible over $K(t)$. We now prove that a similar result without
the completeness assumption on $K$. 

 We begin with the  following

\paragraph*{Lemma 4.4.}
Let $R$ be a  discrete valuation ring and $\pi \in R$ a
parameter. Let $K$ be the field of  fractions of $R$ and $k$ the
residue filed of $R$. Let $F = K(t)$ be the rational function field
in one variable over $K$. Let $n$ be a natural number which is
coprime to char$(k)$. Assume that $K$ contains a primitive
$n^{th}$ root of unity and there exists a $\lambda_{0}\in k^*$  such
that $[k(\sqrt[n]{\lambda_{0}}) : k]=n$. Then $(t,\pi - \lambda
t)_{n}\otimes_F (t+1,\pi)_{m}$ is division over $F$ for any $\lambda
\in R$ which maps to $\lambda_{0}$.

\paragraph*{Proof} Since $\pi$ is a parameter in $R$, the localisation
$R[t]_{(\pi)}$ of $R[t]$ at the prime ideal $(\pi)$ is a discrete
valuation ring. Let $\nu$ be the discrete valuation on $F$ given the
discrete valuation ring $R[t]_{(\pi)}$ and $F_{\nu}$ be the
completion of $F$ at $\nu$. Since the residue field at $\nu$ is 
$k(t)$, we have $[F_{\nu}(\sqrt[m]{t+1}) : F_{\nu} ]=m$.
Let $w$ be the extension of $\nu$ to
$F_{\nu}(\sqrt[m]{t+1})$. To show that $(t,\pi - \lambda
t)_{n}\otimes_F (t+1,\pi)_{m}$ is division, it is enough to show that
$(t,\pi - \lambda t)_{n}\otimes_F (t+1,\pi)_{m} \otimes_F F_{\nu}$ is
division. Since $F_{\nu}$ is complete and  $[F_{\nu}(\sqrt[m]{t+1})
:F_{\nu} ]=m$, by (2.1), it is enough to show that $(t,\pi - \lambda
t)_{n}\otimes F_{\nu}(\sqrt[m]{t+1})$ is division. Since $t$ and $\pi
- \lambda t$ are units at $\nu$, the algebra $(t,\pi - \lambda
t)_{n} \otimes F_{\nu}(\sqrt[m]{t+1})$ is unramified at $w$. Thus it
is enough to show that its image  $ (t, -\lambda_0t)_n = (t, \lambda_0)_n$ is division
over the residue filed $k(t)(\sqrt[m]{t+1})$.   Let $\gamma$ be the
discrete valuation on $k(t)$ given by $t$ and $\tilde{\gamma}$ be
the extension of $\gamma$ to $k(t)(\sqrt[m]{t+1})$. Since $t+1$ is
an $m^{th}$ power in the completion of $k(t)$ at $\gamma$, the
completion of $k(t)(\sqrt[m]{t+1})$ at $\tilde{\gamma}$  is
$k((t))$.  It is enough show that  $(t, \lambda_{0})_{n}$ is
division over $k((t))$. Since $\lambda_{0}\in k^*$ is an element of
order $n$, $(t, \lambda_{0})_{n}$ is division over $k((t))$. \hfill
$\Box$

\paragraph*{Theorem 4.5.}  Let $K$ be a field with a discrete
valuation (not necessarily complete) and $k$ its    
  residue field.  Let $G =  
\mathbb{Z}/l_{1}\mathbb{Z} \times
\mathbb{Z}/l_{2}\mathbb{Z} \times \mathbb{Z}/l_{3}\mathbb{Z} \times
\mathbb{Z}/l_{4}\mathbb{Z}$.
Suppose that the $n = l_1l_2l_3l_4$ 
is coprime to char$(k)$, $K$ contains a primitive $n^{th}$ root of
unity  and there is a $\lambda_0  \in k$ such that 
$[k(\sqrt[l_1l_2]{\lambda_0} : k] = l_1l_2$. 
Then $G$ is admissible over $K(t)$.

\paragraph*{Proof.}    Let $R$ be the ring of integers in $K$
and $\lambda \in R$  mapping to $\lambda_0$ in $k$. Let $\pi \in R$ be
a parameter.   Let  $D_{1}=(t, \pi - \lambda t)_{l_{1}l_{2}}$ and $D_{2}=(t+1,
\pi)_{l_{3}l_{4}}$. Then, by (4.4), $D_1 \otimes D_2$ is a division
algebra over $K(t)$.  Let
 $L_{1}=K(t)(y,z)$ where $y^{l_{1}}=t$,
$z^{l_{2}}= \pi -\lambda t$.   Then $L_{1}$ is a  maximal
subfield of $D_1$.  Let $L_{2}=K(t)(r,s)$ where $y^{l_{3}}=t+1$,
$z^{l_{4}}= \pi$.  Then $L_{2}$ is a  maximal subfield of $D_2$.   
Therefore $L_{1}\otimes _{K(t)} L_{2}$ is a maximal subfield of  
$D_{1}\otimes D_{2}$. Since
$Gal(L_{1}\otimes_{K(t)} L_{2}/K(t))=G$, $G$ is admissible over
$K(t)$. \hfill $\Box$

\paragraph*{Corollary 4.6.}
Let $p$ be a prim and $l_1, l_2, l_3, l_4$ be natural numbers which are coprime
to  $p$. Let $G = \mathbb{Z}/l_{1}\mathbb{Z} \times
\mathbb{Z}/l_{2}\mathbb{Z} \times \mathbb{Z}/l_{3}\mathbb{Z} \times
\mathbb{Z}/l_{4}\mathbb{Z}$.  Let $n = l_1l_2l_3l_4$ and $\zeta$ a
primitive $n^{th}$ root of unity.  Then $G$ is admisible over
$\mathbb{Q}(\zeta)(t)$.

\paragraph*{Remark 4.7.} For $p$, $l_1, l_2, l_3, l_4$, $G$ as in
(4.6),  $G$ is admissible over $\mathbb{Q}_p(t)$ by (4.3). However
(4.5) gives an explicite constriction of a division algebra over
$\mathbb{Q}_p(t)$ and a maximal subfield with Galois group $G$. 

\newpage

\section*{References}

\begin{enumerate}

\item[{[A]}]
\text{Shreeram S. Abhyankar}, \textit{Resolution of singularities of
  algebraic surfaces}.  Algebraic Geometry (Internat.  Colloq.,
Tata Inst. Fund. Res.,  Bombay,1968 ), pp, 1-11, 1969, Oxford Univ. Press, London.

\item[{[AG]}] \text{ Auslander  Maurice and  Goldman Oscar}, 
\textit{The Brauer group of a commutative ring},
Trans. Amer. Math. Soc. {\bf 97} (1960) 367–409.

\item[{[CS]}] \text{David Chillag and Jack Sonn},
  \textit{Sylow-metacyclic groups and $\mathbb{Q}$-admissibility},
  Israel J. Math. {\bf 40} (1981), 307–323.

\item[{[C]}] \text{M. Cipolla},  \textit{Remarks on the lifting of
    algebras over henselian pairs}, 
Mathematische Zeitschrift, {\bf 152} (1977), 253–257.  

\item[{[CTPS]}]
\text{J.-L. Colliot-Thélène, R. Parimala, V. Suresh},
\textit{Patching and local global principles for homogeneous spaces
over function fields of p-adic curves}, to appear in Comm.Math.Helv.

\item[{[CTS]}] \text{ J.-L. Colliot-Th\'el\`ene et J.-J.  Sansuc},  
\textit{Fibr\'es quadratiques et composantes connexes r\'eelles},
Math. Ann.  {\bf 244}  (1979),   105--134. 

\item[{[FS]}] \text{ Burton Fein and Murray Schacher},  \textit{$Q(t)$
    and $Q((t))$-Admissibility of Groups of Odd Order},  Proceedings of
  the AMS,  {\bf 123} (1995),  1639-1645

\item[{[Fe1]}] \text{Walter Feit}, \textit{ The $K$-admissibility of $2A_6$
    and $2A_7$},  Israel J. Math. {\bf 82} (1993),  141–156.  

\item[{[Fe2]}] \text{Walter Feit }, \textit{  $SL(2, 11)$ is
    $\mathbb{Q}$-admissible},  J. Algebra {\bf 257} (2002),   244–248. 

\item[{[Fe3]}] \text{Walter Feit }  \textit{ $PSL_2(11)$ is admissible
    for all number fields},  Algebra, arithmetic and geometry with
  applications (West Lafayette,  IN, 2000),  Springer, Berlin, 2004, pp. 295–299.  

\item[{[FF]}] \text{Paul Feit and Walter Feit}, \textit{ The
    K-admissibility of $SL(2,5)$}, 
Geom. Dedicata {\bf 36} (1990),   1–13.  

\item[{[FV]}] \text{Walter Feit and Paul Vojta}, \textit{Examples of
    some $\mathbb{Q}$-admissible groups},  J. Number Theory {\bf 26} 
(1987),n0.2, 210-226.

\item[{[G]}] \text{Wulf Dieter Geyer}, \textit{An example to Dan's
    talk}, GTEM Summer School,  Geometry and Arithmetic around Galois
  Theory, Galatasaray University,  Istanbul, June 8-19, 2009.

\item[{[GS]}] \text{Philippe Gille and Tammas Szamuely},
  \textit{Central simple algebras and Galois cohomology},  Cambridge
  Studies  in Advanced Mathematics, vol.101, Cambridge University
  Press, Cambridge, 2006.

\item[{[HH]}] \text{David Harbater and Julia Hartmann},
  \textit{Patching over fields},  Israel J.  Math. 
  {\bf 176} (2010), 61-107.

\item[{[HHK1]}] \text{David Harbater, Julia Hartmann and Daniel
    Krashen}, \textit{Applications of patching to quadratic  forms and
    central simple algebras}, Invent.  Math.  {\bf 178}
  (2009),  231-263.  

\item[{[HHK2]}] \text{David Harbater, Julia Hartmann and Daniel Krashen},
\textit{Patching subfields of division algebras}, Trans. A.M.S. {\bf 363} (2011),  3335-3349.

\item[{[HHK3]}] \text{David Harbater, Julia Hartmann and Daniel Krashen},
\textit{Local-global principles for torsors over arithmetic curves}, arXiv:1108.3323v2

\item[{[JW]}] \text{Bill Jacob and Adrian Wadsworth}, \textit{Division
    algebras over Henselian fields},   J. Algebra {\bf 128} (1990),
  126-179. 

\item[{[J]}] \text{Nathan Jacobson}, \textit{Finite-dimensional
    division algebras over fields}, Springer-Verlag, Berlin, 1996

\item[{[KOS]}] \text{M.-A. Knus,  M.  Ojanguren and D. J. Saltman},
    \textit{On Brauer groups in characteristic p},
 Brauer groups (Proc. Conf., Northwestern Univ., Evanston, Ill.,
 1975), pp. 25–49.  Lecture Notes in Math., {\bf 549}, Springer, Berlin, 1976.

\item[{[L]}] \text{Joseph Lipman}, \textit{Introduction to resolution
    of singularities}, in Algebraic geometry, pp,
  187-230. Amer.Math.Soc., Providence,R.I,1975.

\item[{[P]}] \text{Richard S.Pierce}, \textit{Associative algebras},
  Springer-Verlag, New York, 1982. Studies in the History of Modern
  Science.9.

\item[{[S1]}] \text{Davis J. Saltman} \textit{Division algebras over
    p-adic curves},  J. Ramanujan Math. Soc. {\bf 12}, 25–47 (1997)

\item[{[S2]}] \text{David J.Saltman}, \textit{Lectures on division
    algebras}, CBMS Regional Conference Series,no.94,  American
  Mathematical Society, Providence, RI, 1999.

\item[{[Sc]}] \text{Murray M. Schacher}, \textit{Subfields of division
    rings I},   J. Algebra {\bf 9} (1968), 451-419.

\item[{[Sch]}] \text{Scharlau, W}, \textit{Quadratic and Hermitian
    Forms}, Grundlehren der Math. Wiss., 270, Berlin,  Heidelberg, New
  York 1985.

\item[{[Se]}] \text{Jean-Pierre Serre}, \textit{Local fields},
  Springer-Verlag, New York, 1979.

\item[{[So]}] \text{Jack Sonn}, \textit{$\mathbb{Q}$-admissibility of
    solvable groups},  J. Algebra {\bf 84} (1983),   411–419.

\item[{[W]}] \text{Adrian R. Wadsworth}, \textit{Valuation theory on
    finite dimensional division algebras,  in valuation theory and its
    applications}, vol I,  Fields Inst.Commun., vol.32,
  Amer.Math.Soc., Providence, RI,2002, pp.385-449. 

\end{enumerate}

\noindent Department of Mathematics and Statistics \\
University of Hyderabad \\ Gahcibowli
\\ Hyderabad - 500046\\ India \\

and \\

\noindent Department of Mathematics and Computer Science  \\
Emory University   \\ Atlanta, GA 30322  \\ U.S.A.

\end{document}